\documentclass[reqno]{amsart}
\usepackage{enumerate,color}

\newtheorem{theorem}{Theorem}[section]
\newtheorem{lemma}[theorem]{Lemma}

\theoremstyle{definition}
\newtheorem{remark}[theorem]{Remark}

\begin{document}
\title{Note on the number of divisors of reducible quadratic polynomials}
\author{Adrian W. Dudek}
\address{Cronulla NSW 2230}
\email{awdudek@gmail.com}

\author{{\L}ukasz Pa\'nkowski}
\address{Faculty of Mathematics and Computer Science, Adam Mickiewicz University, Umultowska 87, 61-614 Pozna\'{n}, Poland}
\email{lpan@amu.edu.pl}

\author{Victor Scharaschkin}
\address{Department of Mathematics, University of Queensland, St Lucia, QLD 4072, Australia}
\email{vscharaschkin@gmail.com}

\thanks{The second author was partially supported by the Grant no. 2016/23/D/ST1/01149 from the	National Science Centre.}



\begin{abstract}
In a recent paper, Lapkova uses a Tauberian theorem to derive the asymptotic formula for the divisor sum $\sum_{n \leq x} d( n (n+v))$ where $v$ is a fixed integer and $d(n)$ denotes the number of divisors of $n$. We reprove her result by following a suggestion of Hooley, namely investigating the relationship between this sum and the well-known sum $\sum_{n \leq x} d( n ) d (n+v)$. As such, we are able to furnish additional terms in the asymptotic formula.
\end{abstract}

\maketitle

\section{Introduction}

The problem of estimating the average number of divisors of a polynomial was first investigated in the middle of the last century. For example, Erd\"os \cite{Erdos} proved that for every irreducible polynomial $P(n)$ with integer coefficients we have
\[
\sum_{n\leq x} d(P(n))\asymp x\log x,
\]
where $d(n)$ counts positive divisors of $n$. The exact asymptotic formula for the sum $\sum_{n\leq x} d(P(n))$ where $P(n)$ is a polynomial of degree greater than $2$ is still unknown and seems to be a very difficult problem. However, the case of irreducible quadratic polynomials of degree $2$ has been thoroughly investigated and it is known that
\begin{equation}\label{eq:quadratic}
\sum_{n\leq x} d(an^2+bn+c)\sim \lambda x\log x
\end{equation}
for any irreducible polynomial $ax^2+bx+c$ with integer coefficients, where $\lambda$ depends on $a,b,c$. It is an unpublished result (mentioned by R. Bellman in \cite{Bellman}) due to R. Bellman and H. Shapiro, but the first published proof was given by E.\,J. Scourfield in \cite{Scourfield}. For the case $a=1$, $b=0$ the precise dependence of $\lambda$ on $a,b,c$ was described by C.~Hooley in \cite{Hooley}, and for other cases in the series of papers \cite{McKee1}--\cite{McKee3} by J. McKee. In general $\lambda$ depends on the class number of the quadratic field defined by $P$, and hence does not admit a completely elementary description.

R. Bellman \cite{Bellman} also mentioned that there is an unpublished result due to R.~Bellman and H. Shapiro that \eqref{eq:quadratic} holds with $\log x$ replaced by $\log^2 x$ for reducible quadratic polynomials. The first published proof for $a=1$, $b=0$ and $c=-1$ was given by the first author in \cite{Dudek}, who proved that $\sum_{n\leq x} d(n^2-1) \sim \frac{6}{\pi^2}x\log^2x$. His approach is essentially based on the precise description of the function $\rho_{a,b,c}(n)$ (in the case where $a=1$, $b=0$ and $c=-1$) denoting the number of solutions of the congruence $ax^2+bx+c\equiv 0$ mod $n$ in $\mathbb{Z}_n$, and the fact, inspired by the approach suggested by R. Bellman in \cite{Bellman}, that the left hand side of \eqref{eq:quadratic} can be written as
\[
2x\sum_{n\leq x}\frac{\rho_{a,b,c}(n)}{n} + O\left(\sum_{n\leq x}\rho_{a,b,c}(n)\right).
\]
Very recently, this approach was extended by K. Lapkova in \cite{Lapkova} for the polynomial $P(x) = (x-b)(x-c)$ with $b<c$. Since in this case the polynomial is reducible over $\mathbb{Z}$, it is reducible modulo~$p$ for every~$p$ and so $\rho_{a,b,c}(p^n)=2$ for almost all~$p$.  It turns out that in this case the constant $\lambda$ in~\eqref{eq:quadratic} does not depend on the coefficients of a given polynomial, since she proved that
\begin{equation}\label{eq:Lapkova}
\sum_{n\leq x} d((n-b)(n-c)) \sim \frac{6}{\pi^2} x\log^2x.
\end{equation}
Moreover, in \cite{Lapkova}, K. Lapkova extended, by using a different method, the recent result of M. Cipu and T. Trudgian \cite{CipuTrudgian} concerning the case $-b=c=1$, and gave the explicit upper bound for the left hand side of \eqref{eq:Lapkova} with $-b=c=4^s$, $s\geq 0$, where the fastest growing term is exactly $\frac{6}{\pi^2} x\log^2x$ and agrees with \eqref{eq:Lapkova}. The explicit upper bounds for these kind of sums with $b=-c$ are important in searching  for $D(c^2)$-$m$-tuples, namely sets of positive integers $\{a_1,\ldots,a_m\}$ such that $a_ia_j+c^2$ is a perfect square for all $i,j$ with $1\leq i < j \leq m$.

In the paper we give a more precise asymptotic formula for the sum $\sum_{n\leq x} d((n-b)(n-c))$ and prove the following result.
\begin{theorem}\label{th:main}
For every positive integer $v$ and every $\varepsilon>0$ we have
\[
\sum_{n\leq x} d(n(n+v)) = \frac{6}{\pi^2}x\left(\log^2 x + A_1(v)\log x + A_2(v)\right) + O(x^{2/3+\varepsilon})
\]
with
\begin{align*}
A_1(v) &= 4\gamma - 2 - 4\frac{\zeta'}{\zeta}(2) - 2\sum_{e|v}\frac{\Lambda(e)}{e}\\
A_2(v) &= \left(2\gamma - 1 -2\frac{\zeta'}{\zeta}(2)\right)^2+1 - 4\frac{\zeta''}{\zeta}(2)+4\left(\frac{\zeta'}{\zeta}(2)\right)^2\\
&\quad - \left(4\gamma - 2- 4\frac{\zeta'}{\zeta}(2)\right)\sum_{e|v}\frac{\Lambda(e)}{e} + 2\sum_{e|v}\frac{\Lambda(e)\log e}{e} + \sum_{e|v}\frac{\Lambda_2(e)}{e},
\end{align*}
where $\gamma$ denotes the Euler--Mascheroni constant, $\zeta(s)$ denotes the Riemann zeta-function, 
\[
\Lambda_k(n) = \sum_{d|n}\mu(d)\left(\log\frac{n}{d}\right)^k,
\]
and $\Lambda = \Lambda_1$ (resp. $\mu$) is the classical von Mangoldt (resp. M\"obius) multiplicative function.
\end{theorem}

In order to prove the above theorem, we shall follow Hooley's suggestion (see \cite{Hooley}) to find the relation between the sum $\sum_{n\leq x}d(n(n+v))$ and $\sum_{n\leq x}d(n)d(n+v)$. The latter sum is well investigated and the problem of finding its asymptotic behavior is known as the binary additive divisor problem. It was first studied by Ingham in relation to the fourth moment of the Riemann zeta function in \cite{Ingham0}. Subsequently, in \cite{Ingham}, Ingham proved that
\[
\sum_{n\leq x}d(n)d(n+v) = \frac{6}{\pi^2}\sigma_{-1}(v)x\log^2x + O(x\log x),
\]
where $\sigma_\alpha(n) = \sum_{d|n} d^\alpha$. Ingham's result was improved by Estermann, who showed that
\begin{equation}\label{eq:Estermann}
\sum_{n\leq x}d(n)d(n+v) = \frac{6}{\pi^2}\sigma_{-1}(v)x\left(\log^2x+c_1(v)\log x + c_2(v)\right) + O(x^{11/12}\log^3 x),
\end{equation}
where
\begin{align*}
c_1(v) &= 4\gamma - 2 - 4\frac{\zeta'}{\zeta}(2) - 4\frac{\sigma^{(1)}_{-1}}{\sigma_{-1}}(v)\\
c_2(v) &= \left(2\gamma - 1 -2\frac{\zeta'}{\zeta}(2)\right)^2+1 - 4\frac{\zeta''}{\zeta}(2)+4\left(\frac{\zeta'}{\zeta}(2)\right)^2\\
&\quad - 2\left(4\gamma - 2- 4\frac{\zeta'}{\zeta}(2)\right) \frac{\sigma^{(1)}_{-1}}{\sigma_{-1}}(v) + 4\frac{\sigma^{(2)}_{-1}}{\sigma_{-1}}(v)
\end{align*}
and $\sigma_{\alpha}^{(k)}(n) = \sum_{d|n}d^\alpha\log^k d$. This problem was later investigated by many mathematicians, but the best estimate for the error term in \eqref{eq:Estermann} is due to Deshouillers and Iwaniec \cite{DesIwaniec}, who showed that the error term is $O(x^{2/3+\varepsilon})$ for every $\varepsilon>0$. This error term, through the method of this paper, appears in Theorem \ref{th:main}, notably as the proof of this theorem relies essentially on the following crucial lemma.
\begin{lemma}\label{lem:main}
For every $v>0$ we have
\begin{equation}\label{eq:mainlemma}
\sum_{n\leq x}d(n)d(n+v) = \sum_{e|v}\sum_{n\leq \tfrac{x}{e}}d(n(n+\tfrac{v}{e})).
\end{equation}
and, in consequence
\begin{equation}\label{eq:mainlemma_rev}
\sum_{n\leq x}d(n(n+v)) = \sum_{e|v}\mu(e)\sum_{n\leq \tfrac{x}{e}}d(n)d(n+\tfrac{v}{e}).
\end{equation}
\end{lemma}

\begin{remark} Thus --- as noted --- our proof relies on results on the binary additive divisor problem.  Similar problems about the self correlations of arithmetic functions (such as the left hand side of \eqref{eq:mainlemma}) can be very difficult.  For the M\"{o}bius function, Chowla conjectured $\sum_{n \leq x} \mu(n+a_1)\cdots \mu(n+a_t) = o(x)$.  The case $t=1$ is already equivalent to the Prime Number Theorem.  Larger $t$ values are related to the recent M\"{o}bius disjointness conjecture of Sarnak.  For the von Mangoldt function $\sum_{n \leq x} \Lambda(n) \Lambda(n+2) \sim Ax$ is essentially the twin prime conjecture.
\end{remark}

\begin{remark}
Let us note that there is no serious obstacle to make the implied constant in \eqref{eq:Estermann} explicit. Then one can easily obtain the following explicit upper bound
\[
\sum_{n\leq x} d(n(n+v)) \leq \frac{6}{\pi^2}x\left(\log^2 x + A_1(v)\log x + A_2(v)\right) + A_3(v)x^{11/12}\log^3x,
\]
where the constants $A_1$ and $A_2$ are defined as before and the constant $A_3$ can be explicitly computed from Estermann's proof of \eqref{eq:Estermann}. The above inequality improves known upper bounds in the sense that it holds for general polynomials and, what is more important, the first three leading terms agree with our asymptotic formula in Theorem \ref{th:main}, whereas known results (see \cite[Lemma 5.2]{CipuTrudgian}, \cite[Theorem~3 and Corollary~4]{Lapkova}) give an explicit upper bound where only the first leading term agrees with the asymptotic formula.
\end{remark}

\section{The proof of Lemma \ref{lem:main}}

First, notice that for every multiplicative function $f(n)$ we have for all integers $a,b$ that
\[
f(a)f(b) = f(\gcd(a,b))f(\operatorname{lcm}(a,b)).
\]
In order to see this, it suffices to consider the case when $a=p^\alpha$, $b=p^\beta$ and then use the fact that $\{\min(\alpha,\beta),\max(\alpha,\beta)\} = \{\alpha,\beta\}$. 


Now let us assume that $f(n)$ is a multiplicative function satisfying, for every prime $p$ and every positive integer $n$, the following identity
\begin{equation}\label{eq:recursive}
f(p^{n+1}) = f(p)f(p^n) - f(p^{n-1}).
\end{equation}
We now prove by induction on $\alpha$ that for such a multiplicative function as described the following identity holds for all integers $\alpha,\beta$ with $0\leq \beta\leq \alpha$.
\begin{equation}\label{eq:fromInduction}
\sum_{m=0}^{\beta}f(p^{\alpha+\beta-2m}) = f(p^\alpha)f(p^{\beta}).
\end{equation}
If $\alpha=0$, then $\beta=0$ and \eqref{eq:fromInduction} holds trivially. If $\alpha=1$, then $\beta=0$ or $\beta=1$. The former case holds trivially, whereas the latter case needs \eqref{eq:recursive} with $n=2$.

Now, let us assume that \eqref{eq:fromInduction} holds for $\alpha\leq A$ and all non-negative integers $\beta\leq\alpha$, and consider $\alpha=A+1$. If $\beta\leq A-1$, then our assertion is implied by using \eqref{eq:recursive} twice along with our inductive hypothesis. So it remains to consider the cases $(\alpha,\beta) = (A+1,A)$ and $(\alpha, \beta)=(A+1,A+1)$. In the first case, it suffices to write the left hand side of \eqref{eq:fromInduction} as $f(p)+\sum_{m=0}^{A-1}f(p^{2A+1-2m})$ and apply \eqref{eq:recursive}. In the second case, we firstly apply \eqref{eq:recursive} for $n=1$ to write the left hand side of \eqref{eq:fromInduction} as $f(p)^2+\sum_{m=0}^{A-1}f(p^{2A+2-2m})$, and then apply again \eqref{eq:recursive} for $n=2A+1$.

Next, let us observe that for every multiplicative function satisfying \eqref{eq:recursive} we have
\[
f(a)f(b) = \sum_{e|\gcd(a,b)} f\left(\tfrac{ab}{e^2}\right).
\]
Indeed, the above equation holds trivially when $\gcd(a,b)=1$, so let us assume that $\gcd(a,b) = \prod_{j=1}^kp_j^{\alpha_j}$ for some positive integers $\alpha_j$'s, and $ab = q\prod_{j=1}^kp_j^{\beta_j}$ for some $\beta_j\geq 2\alpha_j$ and some integer $q$ coprime with $p_j$'s. Then \eqref{eq:fromInduction} gives us that
\begin{align*}
\sum_{e|\gcd(a,b)}f\left(\tfrac{ab}{e^2}\right) &= f(q)\sum_{\substack{(a_1,\ldots,a_k)\in\mathbb{Z}\\0\leq a_j\leq \alpha_j}}\prod_{j=1}^kf(p_j^{\beta_j-2a_j}) = f(q)\prod_{j=1}^k\sum_{a_j=0}^{\alpha_j}f(p_j^{\beta_j-2a_j})\\
&=f(q)\prod_{j=1}^kf(p_j^{\alpha_j})f(p_j^{\beta_j-\alpha_j})= f(\gcd(a,b))f(\operatorname{lcm}(a,b))\\
& = f(a)f(b).
\end{align*}

Since the multiplicative function $d(n)$ satisfies \eqref{eq:recursive}, one obtains the following lemma.
\begin{lemma}\label{lem:2}
	Let $v$ be a positive integer. Then
	\[
	d(n)d(n+v) = \sum_{e|\gcd(n,v)}d\left(\frac{n(n+v)}{e^2}\right).
	\]
\end{lemma}

Now we are ready to prove \eqref{eq:mainlemma}. Lemma \ref{lem:2} implies \eqref{eq:mainlemma}, since
	\begin{align*}
	\sum_{f|v}\sum_{n\leq \frac{x}{f}}d\left(n\left(n+\frac{v}{f}\right)\right)&=\sum_{f|v}\sum_{\substack{n\leq x\\f|n}}d\left(\frac{n(n+v)}{f^2}\right) \\
	&=\sum_{e|v}\sum_{\substack{n\leq x\\\gcd(n,v)=e}}\sum_{f|e}d\left(\frac{n(n+v)}{f^2}\right)\\
	&=\sum_{e|v}\sum_{\substack{n\leq x\\\gcd(n,v)=e}}d(n)d(n+v)\\
	&=\sum_{n\leq x} d(n)d(n+v).
	\end{align*}

On the other hand, one can easily deduce from \eqref{eq:mainlemma} that
\begin{align*}
\sum_{e|v}\mu(e)\sum_{n\leq \tfrac{x}{e}}d(n)d(n+\tfrac{v}{e})&=\sum_{e|v}\mu(e)\sum_{f|\tfrac{v}{e}}\sum_{n\leq \tfrac{x}{ef}}d(n(n+\tfrac{v}{ef}))\\
&=\sum_{e'|v}\sum_{n\leq \tfrac{x}{e'}}d(n(n+\tfrac{v}{e'}))\sum_{e|e'}\mu(e)\\
&=\sum_{n\leq x} d(n(n+v)).
\end{align*}

\begin{remark}
The crucial property of a multiplicative function $f(n)$ for the above reasoning is to satisfy \eqref{eq:recursive}. In the literature there are many well-known multiplicative functions satisfying the similar identity
\begin{equation}\label{eq:genRec}
f(p^{n+1}) = f(p)f(p^n) - g(p)f(p^{n-1})
\end{equation}
for a suitable completely multiplicative function $g$. Obviously, from our point of view, the case $g\equiv 0$ is not interesting as it implies that $f$ is completely multiplicative, so let us assume that $g\not\equiv 0$. Then, for example, $\sigma_{\alpha}$ satisifes the above identity with $g(p)=p^\alpha$. Moreover, it was noticed by Ramanujan, and proved by Mordell \cite{Mordell}, that Ramanujan's $\tau$ function satisfies this identity with $g(p)=p^{11}$, and more generally, \eqref{eq:genRec} is true for normalized eigenforms of weight $2k$ with $g(p)=p^{2k-1}$. 

Using a similar argument as above, one can easily show that for every multiplicative function $f(n)$ satisfying \eqref{eq:genRec} we have
\[
\sum_{n\leq x} f(n)f(n+v) = \sum_{e|v}g(e)\sum_{n\leq \tfrac{x}{e}}f\left(n(n+\tfrac{v}{e})\right)
\]
and, in consequence, since every non-zero completely multiplicative function $g$ is inverse to $\mu g$ with respect to the Dirichlet convolution, we have 
\[
\sum_{n\leq x} f(n(n+v)) = \sum_{e|v}\mu(e)g(e)\sum_{n\leq \tfrac{x}{e}}f(n)f\left(n+\tfrac{v}{e}\right).
\]
Hence, the asymptotic behavior of $\sum_{n\leq x} f(n(n+v))$ can be deduced from the behavior of $\sum_{n\leq x} f(n)f(n+v)$ and vice-versa. For example, one can easily deduce from \cite{Halberstam} that for $\alpha>0$ we have
\[
\sum_{n\leq x}\sigma_{\alpha}(n(n+v)) = \frac{1}{2\alpha+1}\frac{\zeta(\alpha+1)^2}{\zeta(2\alpha+2)}x^{2\alpha+1}\sum_{d|v}d^{-2\alpha-1}\sum_{e|d}\mu(e)e^\alpha + O(x^\omega\log^c x),
\]
where $\omega=2\alpha+1-\min(\alpha,1)$ and $c=\begin{cases}0,&\alpha>1,\\1,&\alpha<1,\\2,&\alpha=1.\end{cases}$
\end{remark}

\section{The proof of the theorem}

First let us note that
\begin{equation}\label{eq:SigmaLambda}
\sum_{e|v}\frac{\mu(e)}{e}\sigma_{-1}^{(k)}(\tfrac{v}{e}) = \sum_{d|v}\frac{\Lambda_k(d)}{d}
\end{equation}
and
\begin{equation}\label{eq:SumLambda^k}
\sum_{k=0}^n\binom{n}{k}\sum_{e|v}\frac{\mu(e)}{e}\sigma_{-1}^{(k)}(\tfrac{v}{e}) (\log e)^{n-k} = \sum_{d|v}\frac{(\log d)^n}{d}\sum_{e|d}\mu(e)=\begin{cases}1,&n=0,\\0,&n\geq 1.\end{cases}
\end{equation}

Note that \eqref{eq:Estermann} together with \eqref{eq:mainlemma_rev} and \eqref{eq:SumLambda^k} for $n=0$ gives that
\[
\sum_{n\leq x}d(n(n+v))\sim \frac{6}{\pi^2}\sum_{e|v}\mu(e)\sigma_{-1}(\tfrac{v}{e})\frac{x}{e}\log^2 x = \frac{6}{\pi^2}x\log^2 x.
\]
Next, combining \eqref{eq:Estermann} with \eqref{eq:mainlemma_rev} yields
\begin{align*}
A_1(v) &= \frac{6}{\pi^2}\sum_{e|v}\frac{\mu(e)}{e}\sigma_{-1}(\tfrac{v}{e})(c_1(\tfrac{v}{e}) - 2\log e)\\
& = \frac{6}{\pi^2}\left( 4\gamma - 2 - 4\frac{\zeta'}{\zeta}(2) - 4\sum_{e|v}\frac{\mu(e)}{e}\sigma_{-1}^{(1)}(\tfrac{v}{e}) - 2\sum_{e|v}\frac{\mu(e)}{e}\sigma_{-1}(\tfrac{v}{e})\log e\right).
\end{align*}
Thus, \eqref{eq:SigmaLambda} and \eqref{eq:SumLambda^k} for $n=1$ give
\begin{align*}
A_1(v) &= \frac{6}{\pi^2}\left( 4\gamma - 2 - 4\frac{\zeta'}{\zeta}(2) - 2\sum_{e|v}\frac{\mu(e)}{e}\sigma_{-1}^{(1)}(\tfrac{v}{e}) \right)\\
&= \frac{6}{\pi^2}\left( 4\gamma - 2 - 4\frac{\zeta'}{\zeta}(2) - 2\sum_{e|v}\frac{\Lambda(e)}{e} \right).
\end{align*}
Similarly one can compute $A_2(v)$. First let us note that
\begin{align*}
A_2(v) &= \frac{6}{\pi^2}\sum_{e|v}\frac{\mu(e)}{e}\sigma_{-1}(\tfrac{v}{e})\left(\log^2e - c_1(\tfrac{v}{e})\log e + c_2(\tfrac{v}{e})\right)\\
&=\frac{6}{\pi^2}\Bigg(\left(2\gamma - 1 -2\frac{\zeta'}{\zeta}(2)\right)^2+1 - 4\frac{\zeta''}{\zeta}(2)+4\left(\frac{\zeta'}{\zeta}(2)\right)^2\\
&\quad +\sum_{e|v}\frac{\mu(e)}{e}\sigma_{-1}(\tfrac{v}{e})\log^2 e+ 4\sum_{e|v}\frac{\mu(e)}{e}\sigma^{(1)}_{-1}(\tfrac{v}{e})\log e+4\sum_{e|v}\frac{\mu(e)}{e}\sigma^{(2)}_{-1}(\tfrac{v}{e})\\
&\quad - \left(4\gamma - 2- 4\frac{\zeta'}{\zeta}(2)\right) \sum_{e|v}\frac{\mu(e)}{e}\sigma_{-1}(\tfrac{v}{e})\log e\\
&\quad - 2\left(4\gamma - 2- 4\frac{\zeta'}{\zeta}(2)\right) \sum_{e|v}\frac{\mu(e)}{e}\sigma^{(1)}_{-1}(\tfrac{v}{e})\Bigg).
\end{align*}
Then, as in the case of $A_1(v)$, we see that the last two summands give 
\[
-\left(4\gamma - 2- 4\frac{\zeta'}{\zeta}(2)\right)\sum_{e|v}\frac{\Lambda(e)}{e}.
\]
Finally, \eqref{eq:SumLambda^k} for $n=2$ together with \eqref{eq:SigmaLambda} and the fact that
\[
\sum_{e|v}\frac{\mu(e)}{e}\sigma_{-1}(\tfrac{v}{e})\log^2 e+ \sum_{e|v}\frac{\mu(e)}{e}\sigma^{(1)}_{-1}(\tfrac{v}{e})\log e=-\sum_{e|v}\frac{\Lambda(e)\log e}{e}
\]
completes the proof.

\end{document}